\documentclass[pre,twocolumn,showpacs]{revtex4}     
\usepackage{epsfig,amsmath} 

\begin{document}

\title{Super-high-efficiency approximate calculation of series sum and discrete Fourier transform} 

\author{Xin-Zhong Yan}
\affiliation{Institute of Physics, Chinese Academy of Sciences, P.O. Box 603, Beijing 100190, China}
 
\date{\today}
 
\begin{abstract}
We present a super-high-efficiency approximate computing scheme for series sum and discrete Fourier transform. The summation of a series sum or a discrete Fourier transform is approximated by summing over part of the terms multiplied by corresponding weights. The calculation is valid for the function under the transform being piecewise smooth in the continuum variable. The scheme reduces significantly the requirement for computer memory storage and enhances the numerical computation efficiency without losing accuracy. By comparing with the known results of examples, we show the accuracy and the efficiency of the scheme. The efficiency can be higher than $10^6$ for the examples.
\end{abstract}


\maketitle

The discrete Fourier transform (DFT) is defined as
\begin{eqnarray}
F(k) = \sum_{n=n_a}^{n_b}f(n)\exp(-ikn) \label{frt}
\end{eqnarray}
where $f(n)$ as a function of integer $n$ is defined in the range $n_a \le n \le n_b$, and $k$ is a real parameter in the range $(-\pi,\pi)$. For $k = 0$, the calculation reduces to the sum of the series of $f(n)$. The series sum or the DFT is a fundamental calculation appearing widely in problems of science as well as in technology \cite{Strang}. In most problems, $f(n)$ cannot be represented by analytical formula but given by numerical data and we encounter the following difficulties: (1) the range $[n_a, n_b]$ is so large that the memory volume for the data of the function exceeds the storage limit of a computer, (2) even no storage problem, the computation may be very time consuming. Generally, the function $f$ may not only depend on $n$ but also be function of other variables $p_1,p_2,\cdots$, $f = f(n,p_1,p_2,\cdots)$, and need large memory to occupy. Moreover, when the transform is involved in integral equations \cite{Fetter,Baym,Bickers,Pao,Monthoux} that are solved by iterations, the transform needs to be calculated again and again in the iteration process. Therefore, a numerical scheme for the DFT that reduces the storage requirement and accelerates the computation process without losing the accuracy is very desirable.

Here, we present a numerical scheme for the series sum and the DFT. We examine its accuracy and efficiency by comparing the scheme with the known results for some examples. 

To reduce the intensive computation that sums over $n$ term by term in Eq. (\ref{frt}), we intend to get a formula that sums over some selected terms. Since the factor $\exp(-ikn)$ sensitively depends on $n$, we need to treat the summation by approximation with cautious. We suppose $f(x)$ as a function of the continuum variable $x$ is continuous and smooth at least piecewisely. Then, we can take parabolic expansion for $f(x)$ in small regions where the function is smooth. With the expansion, the summation can be carried out analytically. This is the basic idea of our scheme. 

Suppose the function $f(x)$ is smooth in a neighborhood of the number $n_1$, we expand $f(x)$ as
\begin{eqnarray}
f(x) \approx f(n_1) + c_1(x-n_1) + c_2(x-n_1)^2 \nonumber
\end{eqnarray}
where the coefficients $c_1$ and $c_2$ are determined by the function values $f(n_2)$ and $f(n_3)$ at other two points $n_2$ and $n_3$ within this neighborhood. The coefficients are determined as
\begin{eqnarray}
c_1 &=& \frac{n_3-n_1}{n_3-n_2}\frac{f(n_2) -f(n_1)}{n_2-n_1}-\frac{n_2-n_1}{n_3-n_2}\frac{f(n_3) -f(n_1)}{n_3-n_1}
, \nonumber\\
c_2 &=& \frac{f(n_3)-f(n_1)}{(n_3-n_1)(n_3-n_2)}-\frac{f(n_2)-f(n_1)}{(n_2-n_1)(n_3-n_2)}. \nonumber
\end{eqnarray}
At the integer number $n$, we have
\begin{eqnarray}
f(n) \approx f(n_1) + c_1(n-n_1) + c_2(n-n_1)^2. \label{fj}
\end{eqnarray}
We hereafter suppose the three numbers are in order $n_1 < n_2 < n_3$ and are integers.

We now consider the following summation with $n$ running over the small range $[n_1,n_3-1]$, 
\begin{eqnarray}
F(k;n_1,n_3) = \sum_{n=n_1}^{n_3-1}f(n)\exp(-ikn). \label{sum}
\end{eqnarray}
Using the expansion Eq. (\ref{fj}) for $f(n)$, we can get approximately an analytical formula for $F(k;n_1,n_3)$. We need the following summation 
\begin{eqnarray}
S_1(k) &=& \sum_{n=n_1}^{n_3-1}\exp(-ikn) \nonumber\\
&=& \exp(-ikn_1) \frac{1-\exp[-ik(n_3-n_1)]}{1-\exp(-ik)}  \nonumber\\
&\equiv &  \exp(-ikn_1)y(k)           \label{s1} \\
S_2(k) &=& \sum_{n=n_1}^{n_3-1}(n-n_1)\exp(-ikn) \nonumber\\
&=& i\exp(-ikn_1)y'(k) \label{s2}\\
S_3(k) &=& \sum_{n=n_1}^{n_3-1}(n-n_1)^2\exp(-ikn) \nonumber\\
&=& -\exp(-ikn_1)y''(k), \label{s3}
\end{eqnarray}
where the function $y(k)$ and its derivatives are given by
\begin{eqnarray}
y(k) &=& \frac{1-\exp[-ik(n_3-n_1)]}{1-\exp(-ik)},    \nonumber\\
y'(k) &=& -\frac{1-\exp[-ik(n_3-n_1)]}{4\sin^2(k/2)}\nonumber\\
 && ~~~~+i\frac{n_3-n_1}{2\sin(k/2)}\exp[-ik(n_3-n_1-1/2)] \nonumber\\
y''(k) &=& \frac{(n_3-n_1+1/2)\exp[-ik(n_3-n_1)]-1/2}{2\sin^2(k/2)}\nonumber\\
~~~~&&+i\frac{1-\exp[-ik(n_3-n_1)]}{4\sin^3(k/2)}\exp(-ik/2)\nonumber\\
~~~~&&+i\frac{(n_3-n_1)^2}{2\sin(k/2)}\exp[-ik(n_3-n_1-1/2)]. \nonumber
\end{eqnarray}
Substituting the reults into Eq. (\ref{sum}), we get
\begin{eqnarray}
F(k;n_1,n_3) \approx f(n_1)S_1(k)+c_1S_2(k)+c_2S_3(k).     \label{sm1}
\end{eqnarray}
Using the expressions for $c_1$ and $c_2$, we write the above result as
\begin{eqnarray}
F(k;n_1,n_3) &\approx & w_1(k,n_1,n_2,n_3)f(n_1)\exp(-ikn_1)\nonumber\\
  &&+w_2(k,n_1,n_2,n_3)f(n_2)\exp(-ikn_2)\nonumber\\
  &&+w_3(k,n_1,n_2,n_3)f(n_3)\exp(-ikn_3) \nonumber
\end{eqnarray}
where the weight functions $w_{1,2,3}(k,n_1,n_2,n_3)$ are given by
\begin{eqnarray}
w_1 &=& y(k) +\frac{[(2n_1-n_2-n_3)y'(k)+y''(k)]}{(n_2-n_1)(n_3-n_1)}, \nonumber\\
w_2 &=& \frac{[(n_3-n_1)y'(k)-y''(k)]}{(n_2-n_1)(n_3-n_2)}\exp[ik(n_2-n_1)], \nonumber\\
w_3 &=& \frac{[y''(k)-(n_2-n_1)y'(k)]}{(n_3-n_1)(n_3-n_2)}\exp[ik(n_3-n_1)]. \nonumber
\end{eqnarray}
At $k = 0$, the weights are defined as
\begin{eqnarray}
w_1|_{k=0} &=&\frac{n_3-n_1+1}{6(n_2-n_1)}(3n_2-2n_1-n_3+1), \nonumber\\
w_2|_{k=0} &=&\frac{(n_3-n_1)[(n_3-n_1)^2-1]}{6(n_2-n_1)(n_3-n_2)}, \nonumber\\
w_3|_{k=0} &=&\frac{n_3-n_1-1}{6(n_3-n_2)}(n_1-3n_2+2n_3-1). \nonumber
\end{eqnarray}

We now go back to the Fourier transform defined by Eq. (\ref{frt}). By selecting a sequence of odd number of integers, $n_a = n_1 < n_2 < \cdots < n_{2m+1} = n_b$, we apply the above rule and get
\begin{eqnarray}
F(k) &=& \sum_{\ell=1}^{m}F(k;n_{2\ell-1},n_{2\ell+1})+f(n_b)\exp(-ikn_b) \nonumber\\
&=& \sum_{j=1}^{2m+1}W_j(k)f(n_j)\exp(-ikn_j) \label{fsm}
\end{eqnarray}
with
\begin{eqnarray}
W_1(k) &=& w_1(k,n_1,n_2,n_3) \nonumber\\
W_{2\ell-1}(k)&=& w_1(k,n_{2\ell-1},n_{2\ell},n_{2\ell+1})\nonumber\\
&&~~+w_3(k,n_{2\ell-3},n_{2\ell-2},n_{2\ell-1})\nonumber\\
&& ~~~~~~~~~~~~~~~~~~\ell = 2,3,\cdots,m\nonumber\\
W_{2m+1}(k) &=& w_3(k,n_{2m-1},n_{2m},n_{2m+1})+1 \nonumber\\
W_{2\ell}(k)&=& w_2(k,n_{2\ell-1},n_{2\ell},n_{2\ell+1}),
~~\ell = 1,3,\cdots,m. \nonumber
\end{eqnarray}

For some problems, we may need to take only the sine or cosine transform:
\begin{eqnarray}
F_s(k) = \sum_{n=n_a}^{n_b}f(n)\sin(kn), \nonumber\\
F_c(k) = \sum_{n=n_a}^{n_b}f(n)\cos(kn). \nonumber
\end{eqnarray}
Suppose $f(n)$ is a real function. By separating the real and imaginary parts of Eq. (\ref{fsm}),  we obtain
\begin{eqnarray}
F_s(k) &=&\sum_{j=1}^{2m+1}f(n_j)[W^R_j(k)\sin(kn_j)-W^I_j(k)\cos(kn_j)],\nonumber\\\label{sin}\\
F_c(k) &=&\sum_{j=1}^{2m+1}f(n_j)[W^R_j(k)\cos(kn_j)+W^I_j(k)\sin(kn_j)],\nonumber\\\label{cos}
\end{eqnarray}
where $W^R_j$ and $W^I_j$ are the real and imaginary parts of weight function $W_j$, respectively. Since the transforms are linear in $f$, Eqs. (\ref{sin}) and (\ref{cos}) also valid for complex function $f$.

To test the accuracy and the efficiency of the scheme, we compare the numerical computations and the exact results for three examples in following. 

{\it Example} 1. We first consider the sum of the typical series $1/n^p$ for $n = 1, 2, \cdots$, 
\begin{eqnarray}
\zeta(p) &=& \sum_{n=1}^{\infty}\frac{1}{n^p},    \label{rmzf}
\end{eqnarray}
which is known as the Riemann zeta function. By numerical summation, $\zeta(p)$ is calculated by summing the terms up to a cutoff $N$. The error due to the dropped terms is about $O(N^{1-p})$. To reduce this error less than a small quantity $\delta$, $N^{1-p}<\delta$, we must have $N \sim \delta^{1/(1-p)}$. For $\delta = 10^{-4}$ and $p = 1.5$, one needs to sum $N = 10^8$ terms. Here, for using the above scheme, we select the integers $n_j$ for $1 \le j \le M = 151$ as
\begin{eqnarray}
n_j = 
\begin{cases} [q^{j-1}], &{\rm if~} [q^{j-1}] > j \\
j, &{\rm otherwise} 
\end{cases}
\end{eqnarray}
where $q = 1.15$ and the square bracket means the integer part of the number. We name this sequence of the integers $n_j$'s as the $q$ sequence. The sum is calculated by
\begin{eqnarray}
S = \sum_{j=1}^{M}W_j(0)/n^p_j. \label{ssum}
\end{eqnarray}
The numerical results of this sum and the precise values of $\zeta(p)$ for various $p$ are listed in Table I. The error $\Delta=S-\zeta(p)$ seems satisfactory small. Since the convergence of the sum is worse for smaller $p$, larger cutoff is needed for high accuracy result. In the present calculation, since the cutoff $[q^{M-1}]$ is the same for all the parameters, the error is therefore larger for smaller $p$. The cutoff $[q^{M-1}] \approx 1.27\times 10^9$ has reached the limit of the largest integer of our computer.

\begin{table}[t]
\caption{\label{tab}Numerical sum $S$ compared with the Riemann function $\zeta(p)$ for various $p$. The quantity $\Delta=S-\zeta(p)$ represents the numerical error.}
\vspace{-3mm}
\begin{center}
\begin{ruledtabular}
\begin{tabular}{cccc}
$p$  & $\zeta(p)$ &$S$ & $\Delta$ \\
\hline
1.4 & 3.1055 &3.1048 &-0.0007 \\
1.5 & 2.6124 & 2.6122& -0.0002\\
1.6 &2.2858& 2.2857 & -0.0001\\
1.7 & 2.0543 & 2.0542 & -0.0001\\
1.8 & 1.8822 & 1.8822 & 0.0000\\
2 & 1.6449 & 1.6449 & 0.0000\\
\end{tabular}
\end{ruledtabular}
\end{center}
\end{table}

The efficiency $c$ is defined as the ratio between the total number of terms to be summed by definition (the cutoff number) and the actual number of terms in the calculation. Since the cutoff is $[q^{M-1}]$, the efficiency of this calculation is
\begin{eqnarray}
c = [q^{M-1}]/M \approx 8.4\times 10^6.
\end{eqnarray}

{\it Example} 2. We define function $F(x)$ as,
\begin{eqnarray}
F(x) = \frac{1}{a}+\frac{2}{a}\sum_{n=1}^{\infty}\frac{\cos(n\pi x)}{(n\pi/a)^2+1} \label{ex2}
\end{eqnarray}
where $a$ is a parameter. The function to be cosine-Fourier transformed is 
\begin{eqnarray}
f(n) = 1/[(n\pi/a)^2+1].  \nonumber
\end{eqnarray}
For large $n$, $f(n) \sim (a/n\pi)^2$. By setting the cutoff for the summation as $N = {\rm Max}(300,300a/\pi)$, the error stemming from the dropped terms of $n > N$ is about $300^{-2} \approx 1.1\times 10^{-5}$. For selecting the integers over which the sum is carried out, we note that $f(n)$ is flat for $n << a/\pi$. We therefore choose $m_0$ equal-spaced integers in the range $[1,N_0]$ with $N_0 = [4a/\pi]+1$. The total number of the summation is given as $M = 151$. The integer $m_0$ is set as $m_0=[4M/5]$ for $N_0 > [4M/5]$ or $m_0 = N_0$ for $N_0 \le [4M/5]$. The $M-m_0$ integers in the range $(N_0, N)$ is chosen as
\begin{eqnarray}
n_j = 
\begin{cases} n_{j-1}+1, &{\rm for~} j=m_0+1,\cdots,j_0-1 \\
[qn_{j-1}], &{\rm for ~} j = j_0,\cdots,M
\end{cases}
\end{eqnarray}
where $q = (N/N_0)^{1/(M-m_0)}$ and $j_0$ is the minimum integer that $[qn_{j_0-1}] > n_{j_0-1}$. So far, we have determined all the integers for the summation. Finally, by comparing with Eq. (\ref{cos}), we note the wave number is $k = x\pi$. We can then calculate the summation using the above cosine-Fourier transform rule. The results of the present numerical scheme for $a = 1$, 5, and $10^5$ are shown as the symbols in Fig. 1. On the other hand, the sum can be exactly obtained. The result for the function is 
\begin{eqnarray}
F(x) = \frac{\cosh[a(1-x)]}{\sinh a}, \nonumber
\end{eqnarray}
for $0 < x < 2$. $F(x)$ is a periodic function of $x$ with periodicity 2. The exact function $F(x)$ is shown in Fig. 1 as the solid lines. From Fig. 1, we see that the numerical results are surprisingly good. For the large parameter $a = 10^5$, the cutoff number of the summation needs to be very large to get the sum converged. The present scheme sums only $M = 151$ terms and is very efficient. 

\begin{figure}[t] 
\centerline{\epsfig{file=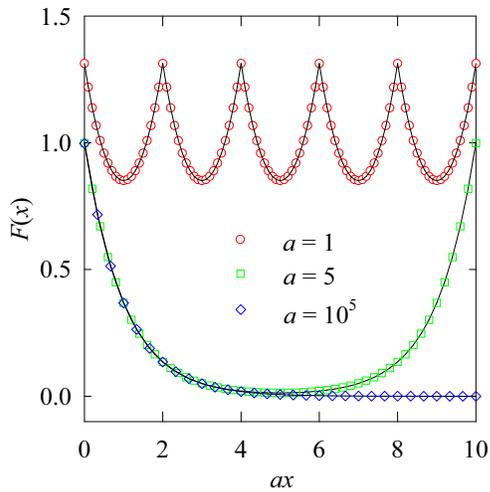,height=7.cm}}
\caption{(color online) $F(x)$ as function of $x$ at parameters $a$ = 1, 5 and $10^5$. Solid lines represent the exact function. Symbols are the numerical results obtained using the present numerical scheme.} 
\end{figure} 

{\it Example} 3. We consider the function $F(x)$ defined as,
\begin{eqnarray}
F(x) = \frac{1}{a}-\frac{2}{a}\sum_{n=1}^{\infty}\frac{\cos(n\pi x)}{(n\pi/a)^2-1} \label{ex3}
\end{eqnarray}
where $a$ is a parameter. The function to be cosine-Fourier transformed
\begin{eqnarray}
f(n) = 1/[(n\pi/a)^2-1]  \nonumber
\end{eqnarray}
is slightly different from the one in Example 2. But $f(x)$ has a singularity at $x = a/\pi$. We therefore need to take cautious in applying the above scheme. Denote $[a/\pi] = N_0$. We separate the summation into two parts respectively over the regions $[1, N_0]$ and $[N_0+1, N]$ with $N = {\rm Max}(1, N_0)\times 10^4$ as the cutoff. Now in each region, since the function is smooth, we can apply the above scheme. For small $N_0$, the calculation in the first region $[1, N_0]$ can be carried out by summing term by term. While for large $N_0$, we choose odd number of integers in this region. Since the function $f(n)$ varies rapidly close to $N_0$, the selected integers should be densely distributed around it. We choose integers using the $q$ sequence with $q = 1.1$ from $N_0$ in the first region. Concretely, we have 
\begin{eqnarray}
n_1 &=& 1   \nonumber\\
n_j &=& N_0+1-\tilde n^{(1)}_{m_1+1-j}, ~~ j= 2,\cdots,m_1  \nonumber
\end{eqnarray}
where $\tilde n^{(1)}_j$'s are the integers determined by the $q=1.1$ sequence, and $m_1$ is an odd number. From the minimum integer $j_0$ that $\tilde n^{(1)}_{j_0} > N_0$, $m_1$ is determined as $m_1 = j_0$ if $j_0$ is odd, or $j_0-1$ if $j_0$ is even. In the second region $[N_0+1, N]$, we choose $M = 151$ integers according to another $q$ sequence,
\begin{eqnarray}
n_j &=& N_0+\tilde n^{(2)}_j, ~~ j= 1,\cdots,M  \nonumber
\end{eqnarray}
where $\tilde n^{(2)}_j$'s are the integers of the sequence of $q = (N-N_0)^{1/(M-1)}$. With the integers so selected, we do the summations according to our scheme in the two regions. The calculated results for $a = \pi/2, 3\pi/2$, and $10^5\pi+\pi/2$ are shown in Fig. 2 as the symbols. The calculation is compared with the exact function
\begin{eqnarray}
F(x) = \sin(ax)-\frac{\cos(ax)}{\tan a}, {\rm ~~~~for~} 0 < x < 2 \nonumber
\end{eqnarray}
which is depicted as the solid lines in Fig. 2. Obviously, the numerical summation by our scheme can very good reproduce the exact result. 

\begin{figure}[t] 
\centerline{\epsfig{file=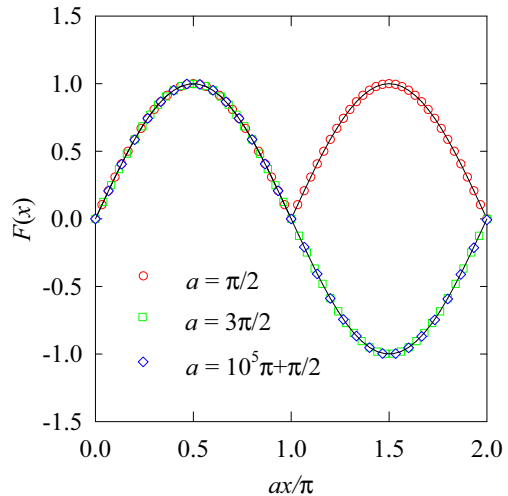,height=7.cm}}
\caption{(color online) $F(x)$ as function of $x$ at parameters $a = \pi/2, 3\pi/2$, and $10^5\pi+\pi/2$. Solid lines represent the exact function. Symbols are the numerical results obtained using the present numerical scheme.} 
\end{figure} 

In the present scheme, the selected integers are in ascending order, but not required to be necessarily equal-space distributed. When they are equal-space selected, the formulas of the present scheme reduce to that of the previous work \cite{Yan1}. In the previous work, we have used an integer-selection plan \cite{Yan2} by which the whole domain for the summation is separated into several subdomains; in each subdomain, the selected integers are equal-space distributed. With such an integer-selection plan, we have obtained compact forms for the weight functions and the Fourier transform. However, the equal-space-integer selection is not convenient for dealing with various problems. For some problems, only finite but very large number of terms is under considerations \cite{Yan3}. This finite number acts as the cutoff. It may not match with the one given by the equal-space-integer selection. Moreover, for a function that may be discontinuous or singular at some points (such as appeared in Example 3), we need to separate the whole domain of the function's definition into several regions where the function is smooth and choose integers densely close to the singular points. For dealing with these problems and for general uses, the present scheme with non-equal-space distributed integers is more convenient.

The DFT can be used to efficiently solve the integral equations in quantum many-particle problems by the Green's function theory \cite{Fetter,Baym,Bickers,Pao,Monthoux}, the quantum Monte Carlo \cite{Foulkes}, and the dynamical mean-field theory \cite{Kotliar}. In quantum mechanics, one encounters the DFT of functions of the quantized energy level or discrete Matsubara frequency \cite{Fetter}, or functions defined in lattices \cite{Hubbard}. It is also useful in many practical applications such as in digital signal processing \cite{Smith}, in image processing \cite{Solomon}, etc. The present scheme may be applicable to the related problems.

By conclusion, we have developed the super-high-efficiency numerical scheme for the series sum and the DFT. The function under transforming is defined at the integers. The function is supposed to be piecewise smooth when its definition is extended to the continuum variable. By the present scheme, the series sum and the DFT are calculated by summing over a number of selected terms with corresponding weights. The corresponding selected integers can be non-equal-space distributed. The integers can be selected densely in a region where the function varies rapidly or sparsely where the function is flat. The scheme reduces significantly the requirement for computer memory storage and thereby enhances the numerical computation efficiency. We have checked the accuracy with several examples. The efficiency can be higher than $10^6$. 

This work was supported by the National Basic Research 973 Program of China under Grants No. 2011CB932702 and No. 2012CB932302.

\end{document}